\documentclass[preprint]{elsarticle}

\usepackage{amsmath}
\usepackage{amsfonts}
\usepackage{amssymb}
\usepackage{amsthm}
\usepackage{bm}
\newtheorem{definition}{Definition}
\newtheorem{proposition}{Proposition}
\newproof{pf}{Proof}
\theoremstyle{plain}
\usepackage{algorithm2e}

\usepackage{graphicx}
\usepackage{epstopdf}

\usepackage{xspace}
\usepackage[font={small,it},labelsep=period,tableposition=top]{caption}
\usepackage{wrapfig}

\usepackage[font={small,it},labelsep=period]{caption}
\usepackage{subcaption}
\usepackage[tableposition=top]{caption}
\usepackage{float}

\usepackage{mathtools}

\usepackage{booktabs}
\usepackage{array}
\usepackage{hyperref}

\begin{document}

\begin{frontmatter}

\title{Analytical solution for heat conduction due to a moving Gaussian heat flux with piecewise constant parameters}

\author[mymainaddress,mysecondaryaddress]{Robert Forslund\corref{mycorrespondingauthor}}
\cortext[mycorrespondingauthor]{Corresponding author}
\ead{robfor@chalmers.se}

\author[mymainaddress]{Anders Snis}

\author[mysecondaryaddress]{Stig Larsson}

\address[mymainaddress]{Arcam EBM, Kroksl{\"a}tts fabriker 27A, SE-431 37 M{\"o}lndal, Sweden}
\address[mysecondaryaddress]{Mathematical Sciences, Chalmers University of Technology and University of Gothenburg, SE-412 96 Gothenburg, Sweden}

\begin{abstract}
We provide an analytical solution of the heat equation in the half-space subject to a moving Gaussian heat flux with piecewise constant parameters. The solution is of interest in powder bed fusion applications where these parameters can be used to control the conduction of heat due to a scanning beam of concentrated energy. The analytical solution is written in a dimensionless form as a sum of integrals over (dimensionless) time. For the numerical computation of these integrals we suggest a quadrature scheme that utilizes pre-calculated look-up tables for the required quadrature orders. Such a scheme is efficient because the required quadrature orders are strongly dependent on the parameters in the heat flux. The possibilities of using the obtained computational technique for the control and optimization of powder bed fusion processes are discussed.
\end{abstract}


\end{frontmatter}


\section{Introduction}
\label{sec1}
Powder bed fusion (PBF) is a particular type of additive manufacturing (AM) where metal powder is melted by a laser or electron beam in a layer-wise fashion to enable the production of geometrically complex parts. AM is a technology that undergoes continuous development towards an industry that is stable and efficient, but there are still issues in terms of processability and quality \cite{marklandkorner, kingetal}. Mathematical models of different levels of complexity are used to gain understanding of the process and to facilitate optimization of the melting. As such, these models provide valuable information and they can complement or replace many types of trial-and-error investigations. 

PBF is a multi-scale process, but macroscopic models are still able to accurately describe the thermal and mechanical behavior of the melting process \cite{riedlbauer}\cite{galati}\cite{jamshidinia}. The essence of these models is that the powder bed is treated as a continuum \cite{marklandkorner} rather than a collection of powder particles. This allows for simulation on large length scales. In the area of continuum modeling of PBF, most contributions are based on the finite element method since it allows for straightforward implementation of different physical phenomena such as cooling, phase transitions, temperature dependencies, among other \cite{schoinochoritisetal}. 

In contrast to finite element methods, analytical solutions of the heat equation are beneficial if the computational effort required to compute them is low. A low computational effort is particularly useful for control and optimization of the melting process, as they rely heavily on an efficient solver. The drawback of analytical solutions is the restrictive assumptions that are required. However, even if detailed physics of the melting process is left out, such solutions can still be used for the purpose of process optimization and verification.

The analytical temperature distribution due to a moving heat source has been derived for a vast variety of problems in both semi-infinite \cite{carslaw, rosenthal, sanders, modest} and finite \cite{carslaw, araya, kidawa, salimi}  three-dimensional domains, and with different heat source models. In this paper the problem under consideration is the heat equation in the lower half-space with a Gaussian distributed heat flux that travels over the surface.  In \cite{clineandanthony} a solution of this problem is obtained via superposition of a series of point heat source solutions in the case of a continuously scanning heat source moving along a single line. These solutions are then expressed in a dimensionless form. It has been shown \cite{sanders} that this solution is a special case of a more general solution and is a valid approximation in the case of strongly absorbing materials where the beam energy is absorbed near the surface. This approximation suits both laser and electron beam applications as long as penetration depths are small.

A later contribution in \cite{john} presents a solution of the heat equation on the half-space for a general heat flux. The solution is obtained via a reformulation of the problem and the use of Duhamel's principle. 
This result is then applied to the traveling Gaussian heat flux problem for a single line with constant beam parameters. The obtained solution coincides with the one presented in \cite{eagerandtsai} and is expressed in a dimensionless form. 

In the present study, we extend this solution to the more general case when the beam path is piecewise straight and beam parameters are piecewise constant in time. In actual PBF processes, the beam parameters can be used to control the temperature during melting and they are often defined in this piecewise constant fashion due to discretization. Therefore, the solution presented here is of particular interest for such applications \cite{snis2}. The solution consists of a sum of integrals over (dimensionless) time. The quadrature orders required to accurately approximate these integrals vary significantly depending on time, beam spot size, and speed. Therefore, the generation of look-up tables that store the required quadrature orders is suggested as a method of speeding up computations. 

The remainder of this paper is organized as follows. In Section \ref{sec2}, we formulate the problem under consideration and derive an analytical solution that is put in dimensionless form. Section \ref{sec3} deals with the numerical computation of said solution. Section \ref{sec4} presents three examples that illustrate the validity of our solution and make use of the quadrature scheme presented in the section prior. Finally, Section \ref{sec5} discusses our solution in the context of optimization of the melting process.


\section{Formulation of the problem and its solution}
\label{sec2}
Consider the heat equation on the lower half space $\Omega = \mathbb{R}^2 \times \mathbb{R}_-$ during a time span $ \mathcal{T} = (0, T]$. Let $\Gamma$ denote the surface boundary $z = 0$ that is subject to a heat flux $\Phi$. Let $u_\mathrm{init}$ denote an initial temperature. With $\hat{z}$ the outward unit normal of $\Omega$, the problem can be written as 

\begin{equation}
\begin{aligned}
\rho c_p \frac{\partial u}{\partial t} - \nabla \cdot ( \lambda \nabla u) &= 0 & \mathrm{in}&  \,\,\,\,\Omega \times \mathcal{T},\\
(\lambda \nabla u) \cdot \hat{z} &= \Phi & \mathrm{on}& \,\,\,\, \Gamma \times \mathcal{T},\\
u(\cdot, 0) &= u_\mathrm{init} & \mathrm{in} & \,\,\,\, \Omega,
\end{aligned}
\label{pde}
\end{equation}

\noindent where $\rho$, $c_p$, $\lambda$ denote density, heat capacity, and thermal conductivity, respectively. These material parameters are assumed to be constant. This gives us a thermal diffusivity $\kappa = \lambda/\rho c_p$. With $\mathbf{x} = (x, y, z)$, the solution of \eqref{pde} is a temperature distribution $u = u(\mathbf{x}, t)$.

From \cite{john}, the solution of \eqref{pde} for a general initial temperature and a general heat flux is given by 
\begin{equation}
u(\mathbf{x}, t) = u^\mathrm{I}(\mathbf{x}, t) + u^\Phi(\mathbf{x}, t), 
\label{sol}
\end{equation}
where
\begin{align}
u^\mathrm{I}(\mathbf{x}, t) &= \left(G( \cdot, t) * u_\mathrm{init}\right)(\mathbf{x}),\label{u11}\\
u^\Phi(\mathbf{x}, t) &= 2 \int_0^t \int_{\mathbb{R}^2} G(x - \xi, y - \eta, z, t-s) \Phi(\xi, \eta, s) \, \mathrm{d}\xi\mathrm{d}\eta \,\mathrm{d}s. \notag
\end{align}
\noindent Here
\[G(\mathbf{x}, t) = \Big( \tfrac{1}{4 \pi \kappa t}\Big)^{3/2} \exp\Big(\scalebox{0.75}[1.0]{\( - \)}\tfrac{|\mathbf{x}|^2}{4\kappa t}\Big)\]
\noindent is the Green's function for three-dimensional diffusion. All convolutions in this paper operate on $\mathbb{R}^3$, 
\begin{align*}
 (f * g)(\mathbf{x}) = \int_{\mathbb{R}^3} f(\bm{\xi})g(\mathbf{x} - \bm{\xi}) \,\mathrm{d}\bm{\xi},
\end{align*}
\noindent with $\bm{\xi} = (\xi, \eta, \zeta)$.

PBF involves a concentrated application of heat via a beam. The beam travels on the surface $\Gamma$ along a path $\mathcal{C}$. The energy from the beam heats the bed due to absorption. The amount of absorbed beam energy depends on several factors such as beam equipment (laser or electron beam), melting conditions (conductive or key hole melting), and material. Here we assume that keyhole effects and penetration depths of the beam are shallow and use a flux to model the contribution of heat. The flux $\Phi$ is modeled by a Gaussian function,
\[\Phi = \Phi(x, y, t)  = \frac{P(t)}{2\pi\sigma(t)^2}\,\mathrm{exp} \Big(\scalebox{0.75}[1.0]{\( - \)}\tfrac{\left(x-x_c\left(t\right)\right)^2+\left(y-y_c\left(t\right)\right)^2}{2\sigma(t)^2}\Big),\]

\noindent where ${\bf{x}}_c = (x_c, y_c, 0) \in \mathcal{C}$ is the position of the center of the beam. The three beam parameters are the absorbed beam power $P(t)$, the effective  beam spot size $\sigma(t)$, and the beam speed $v(t) = |\mathbf{v}(t)|$. The absorbed beam power and the effective beam spot size are regarded as parameters that can be adjusted to mimic a more realistic absorption (see \cite{galati} for a discussion about electron beam heat sources, effective beam diameters and beam efficiency). We have $\mathbf{v} = (v_x, v_y, 0) = (v\cos{\alpha}, v\sin{\alpha}, 0)$, where $\alpha$ is the angle between the $x$-axis and the direction of the path. This angle is known for any $t$ since the beam path $\mathcal{C}$ is pre-set and thus a rotation of axis shows that $v$ uniquely defines the vector $(v_x, v_y, 0)$. The position of the beam ${\bf{x}}_c$ at a given time depends on the speed with which the beam has traveled the path $\mathcal{C}$ up to this time. By using this flux, we make the assumption that surface radiative and evaporative effects are negligible compared to the heat conduction. 

We henceforth also assume that the initial temperature $u_\mathrm{init}$ is constant. This assumption simplifies the exposition as it implies that $\left(G(\cdot, t) * u_\mathrm{init}\right)(\mathbf{x}) = u_\mathrm{init}$. However, it should be noted that a constant initial temperature is not a requirement for the proposed method. A bulk time dependent temperature distribution that for instance mimics the temperature history from previous melted layers can easily be included by summation. 

In the case of constant beam parameters and single line melting, $u^\mathrm{I}(\mathbf{x}, t) = u_\mathrm{init}$ and ${\bf{x}}_c(t)= (v_x t, v_y t, 0)$. Then an explicit evaluation of the spatial integral in the second term of \eqref{sol} gives \cite{john}
\begin{equation}\label{us_single} 
\begin{aligned}
u(\mathbf{x}, t) &= u^\mathrm{I}(\mathbf{x}, t) + u^\Phi(\mathbf{x}, t) \\
& =  u_\mathrm{init} + \int_{\scriptscriptstyle 0}^{\scriptscriptstyle t}{\frac{P\exp\left(\scalebox{0.75}[1.0]{\( - \)}\frac{(x-v_xs)^2+(y-v_ys)^2}{2\sigma^2 + 4\kappa(t-s)}\scalebox{0.75}[1.0]{\( - \)} \frac{z^2}{4\kappa(t-s)}\right)}{(\pi^3\kappa)^{1/2}\rho c_p (t-s)^{1/2}\left(\sigma^2 + 2\kappa(t-s)\right)} \, \mathrm{d} s.} 
\end{aligned}
\end{equation}

\noindent PBF is a dynamic process and to assume constant beam parameters is not a feasible approach since these parameters can, and should, be used to control the melting process. Therefore, we now extend the above solution to the case when beam parameters are piecewise constant and the beam path is not restricted to a single line. Parameter data are often restricted to vary in this discontinuous fashion during melting. The following definition makes the concept of piecewise constant beam parameters more precise.

\begin{definition}\label{def1}
Given times $0 = t_0 < t_1 < \hdots < t_N = T$, we define a partition of $\mathcal{T}$ consisting of $N$ segments
\[(t_{n-1}, t_n] = \big(t_n^\mathrm{i}, t_n^\mathrm{f}\big], \,\,\,\, n = 1, 2, \hdots, N.\]
Index $n$ indicates the $n^{th}$ segment in the partition, and a segment in turn is a collection of the following data:
\begin{itemize}
\item[--] $t_n^\mathrm{i}, \,\,t_n^\mathrm{f}$; an initial time and final time, respectively,
\item[--] $(P_n, \sigma_n, v_n)$; a constant triplet of beam power, spot size, and speed such that 
\[\left(P(t), \sigma(t), v(t)\right) = (P_n, \sigma_n, v_n)    
\,\,\,\,\text{if}\,\,\,\, t \in \big(t_n^\mathrm{i}, t_n^\mathrm{f}\big],
\]
\item[--] $\ell_n$; a line traversed by the beam between times $t_n^\mathrm{i}$ and $t_n^\mathrm{f}$, 
\item[--] $\mathbf{x}_n^\mathrm{i}, \,\,\mathbf{x}_n^\mathrm{f}$; the initial position and final position of $\ell_n$, respectively,
\item[--] $u_n^\mathrm{I}$, $u_n^\Phi$, $u_n$; restrictions of $u^\mathrm{I}$, $u^\Phi$, and $u$ onto segment $n$, respectively, meaning that
\[
\begin{cases}
u(\cdot, t) \,\,= u_n(\cdot, t)\\
u^\mathrm{I}(\cdot, t)  = u_n^\mathrm{I}(\cdot, t)\\
u^\Phi(\cdot, t) = u_n^\Phi(\cdot, t)\\
\end{cases}  \,\text{if}\,\,\,\, t \in \big(t_n^\mathrm{i}, t_n^\mathrm{f}\big].\]
\end{itemize}
\end{definition}

\noindent Note that $t_n^\mathrm{f} = t_{n+1}^\mathrm{i}$.

We now utilize the constant parameters solution \eqref{us_single} in order to express a solution of \eqref{pde} for the practically relevant case when the beam parameters are piecewise constant in the sense of Definition \ref{def1}.
This is done in two steps where we deal with $u^\mathrm{I}$ and $u^\Phi$ separately. Given $t \in \big(t_n^\mathrm{i}, t_n^\mathrm{f}\big]$, the flux term $u^\Phi$ is obtained by translating in space and time. More precisely, with the integrand
\[h(\mathbf{x}, t, s; P, \sigma, v) = \frac{P\exp\left(\scalebox{0.75}[1.0]{\( - \)}\frac{(x-v_xs)^2+(y-v_ys)^2}{2\sigma^2 + 4\kappa(t-s)} \scalebox{0.75}[1.0]{\( - \)} \frac{z^2}{4\kappa(t-s)}\right)}{(\pi^3\kappa)^{1/2} \rho c_p (t-s)^{1/2}\left(\sigma^2 + 2\kappa(t-s)\right)}\]
\noindent and \eqref{us_single}, the flux term $u_n^\Phi$ on segment $n$ starting in $\mathbf{x}_n^\mathrm{i}$ at time $t_n^\mathrm{i}$ can be expressed as
\begin{equation} u_n^\Phi(\mathbf{x}, t) = \int_0^{t-t_n^\mathrm{i}} h(x-x_n^\mathrm{i}, y - y_n^\mathrm{i}, z, t-t_n^\mathrm{i}, s; P_n, \sigma_n, v_n) \,\mathrm{d}s.\label{uns}\end{equation}

\noindent Known substitutions \cite{john} put us in a frame that moves with the beam and let us compute the right-hand side of \eqref{uns} as

\begin{align*}
u_n^\Phi(\bar{\mathbf{x}}_n, \bar{t}_n) & = T_n \int_0^{\bar{t}_n} \tfrac{1}{1+{\bar{s}}^2} \exp\Big(\scalebox{0.75}[1.0]{\( - \)}\tfrac{({\bar{x}_n}+\bar{v}_{n, x}{\bar{s}}^2)^2+{(\bar{y}_n}+\bar{v}_{n, y}{\bar{s}}^2)^2}{1+{\bar{s}}^2} \scalebox{0.75}[1.0]{\( - \)} \tfrac{{\bar{z}_n}^2}{{\bar{s}}^2}\Big) \,\mathrm{d}\bar{s}, \\
\stepcounter{equation}\tag{\theequation}\label{ustrans}
\end{align*}
where 
\begin{align*}
T_n &= P_n/\sqrt{2}\pi^{3/2}\lambda\sigma_n,\,\,\,\,\,\,\,\,\,\,\,\,  \mathrm{[K]} \\
\bar{v}_{n, x} &= v_{n, x}\sigma_n/2\sqrt{2}\kappa, \\
\bar{v}_{n, y} &= v_{n, y}\sigma_n/2\sqrt{2}\kappa, \\
\bar{s} &= \sqrt{2\kappa(t - t_n^\mathrm{i} -s)}/\sigma_n,\\[3pt]
\bar{x}_n &= \big(x - x_n^\mathrm{i} - v_{n, x}(t-t_n^\mathrm{i})\big)/\sqrt{2}\sigma_n, \\
\bar{y}_n &= \big(y - y_n^\mathrm{i} - v_{n, y}(t-t_n^\mathrm{i})\big)/\sqrt{2}\sigma_n, \\
\bar{z}_n &= z/\sqrt{2}\sigma_n, \\
\bar{t}_n &= \sqrt{2\kappa(t-t_n^\mathrm{i})}/\sigma_n.
\end{align*}

\noindent  An overline indicates that a variable is dimensionless. This dimensionless form of the integral is easier to compute numerically.

The initial value term $u_n^\mathrm{I}$ collects the heat contribution due to $u_\mathrm{init}$ as well as the traversing of the beam over all previous segments (with index less than $n$). As a consequence, it is given by the following sum of convolutions,
\begin{equation} u_n^\mathrm{I}(\mathbf{x}, t) = u_\mathrm{init} + \sum_{k=1}^{n-1} \big(G(\cdot, t- t_k^\mathrm{f}) * u_{k}^\Phi(\cdot, t_k^\mathrm{f})\big)(\mathbf{x}). \label{ui}\end{equation}
\noindent We refer to Appendix \ref{app} for a proof of this result. In order to express $u_n^\mathrm{I}$ without spatial integrals we need to compute convolutions on the form
\[ u_{n,k}^\mathrm{I}(\mathbf{x}, t) = \big(G(\cdot, t- t_k^\mathrm{f}) * u_{k}^\Phi(\cdot, t_k^\mathrm{f})\big)(\mathbf{x}) \,\, , \,\,\,\,\, 1 \leq k \leq n-1.\]
\noindent To this end, begin by changing the order of integration,

\begin{align*}
 u_{n,k}^\mathrm{I}(\mathbf{x}, t) &= \!\begin{multlined}[t]\int_{\mathbb{R}^3}
\Big( \tfrac{1}{4 \pi \kappa(t-t_k^\mathrm{f})}\Big)^{3/2} 
\exp\Big(\scalebox{0.75}[1.0]{\( - \)}\tfrac{|\mathbf{x}-\bm{\xi}|^2}{4\kappa(t-t_k^\mathrm{f})}\Big) \\ \cdot\left(\int_0^{t_k^{f}-t_{k}^\mathrm{i}} h(\xi - x_{k}^\mathrm{i}, \eta - y_{k}^\mathrm{i}, \zeta, t_k^\mathrm{f} - t_{k}^\mathrm{i}, r; P_{k}, \sigma_{k}, v_{k}) \,\mathrm{d}r \right)\, \mathrm{d}\bm{\xi}\end{multlined} \\
&= \int_{0}^{t_k^{f}-t_{k}^\mathrm{i}}\int_{\mathbb{R}^3}
\tfrac{h(\xi - x_{k}^\mathrm{i}, \eta - y_{k}^\mathrm{i}, \zeta, t_k^\mathrm{f} - t_{k}^\mathrm{i}, r; P_{k}, \sigma_{k}, v_{k})}{\left(4 \pi \kappa(t-t_k^\mathrm{f})\right)^{3/2}} 
\exp\Big(\scalebox{0.75}[1.0]{\( - \)}\tfrac{|\mathbf{x}-\bm{\xi}|^2}{4\kappa(t-t_k^\mathrm{f})}\Big)\, \mathrm{d}\bm{\xi} \,\mathrm{d}r.
\end{align*}

\noindent A change of variables gives
\[
u_{n,k}^\mathrm{I}(\bar{\mathbf{x}}_k, \bar{t}_k) = \frac{T_k}{\pi^{3/2}\,\bar{t}_k^3} \int_0^{\bar{t}_k^\mathrm{f}} \frac{1}{1 + \bar{r}^2}J_{\bar{\xi}}(\bar{{x}}_k, \bar{t}_k, \bar{r})\,J_{\bar{\eta}}(\bar{{y}}_k, \bar{t}_k, \bar{r})\,J_{\bar{\zeta}}(\bar{{z}}_k, \bar{t}_k, \bar{r}) \,\,\mathrm{d}\bar{r}, 
\]
\noindent where $J_{\bar{\xi}}$, $J_{\bar{\eta}}$, and $J_{\bar{\zeta}}$ are integrals over space that can be evaluated as
\begin{align*}
J_{\bar{\xi}}(\bar{{x}}_k, \bar{t}_k, \bar{r}) 
&= \int_{-\infty}^{\infty} \exp \Big( \scalebox{0.75}[1.0]{\( - \)} \tfrac{(\bar{x}_k-\bar{\sigma}_k\bar{\xi})^2}{\bar{t}_k^2} \scalebox{0.75}[1.0]{\( - \)} \tfrac{(\bar{\xi} + \bar{v}_{k, x}\bar{r}_k^2)^2}{1+\bar{r}_k^2}\Big)\,\mathrm{d}\bar{\xi} \\
& = \sqrt{\frac{\pi\, \bar{t}_k^2(1+\bar{r}^2)}{1 + \bar{r}^2 + \bar{t}_k^2}}\exp\Big(\scalebox{0.75}[1.0]{\( - \)}\tfrac{(\bar{x}_k + \bar{v}_{k, x}\bar{r}^2)^2}{1 + \bar{r}^2 + \bar{t}_k^2})\Big), \\
\intertext{}
J_{\bar{\eta}}(\bar{{y}}_k, \bar{t}_k, \bar{r}) 
&= \int_{-\infty}^{\infty} \exp \Big( \scalebox{0.75}[1.0]{\( - \)} \tfrac{(\bar{y}_k-\bar{\sigma}_k\bar{\eta})^2}{\bar{t}_k^2} \scalebox{0.75}[1.0]{\( - \)} \tfrac{(\bar{\eta} + \bar{v}_{k, y}\bar{r}_k^2)^2}{1+\bar{r}_k^2}\Big)\,\mathrm{d}\bar{\eta} \\
& = \sqrt{\frac{\pi\, \bar{t}_k^2(1+\bar{r}^2)}{1 + \bar{r}^2 + \bar{t}_k^2}}\exp\Big(\scalebox{0.75}[1.0]{\( - \)}\tfrac{(\bar{y}_k + \bar{v}_{k, y}\bar{r}^2)^2}{1 + \bar{r}^2 + \bar{t}_k^2}\Big), \\
J_{\bar{\zeta}}(\bar{{z}}_k, \bar{t}_k, \bar{r}) 
&= \int_{-\infty}^{\infty} \exp \Big( \scalebox{0.75}[1.0]{\( - \)} \tfrac{(\bar{z}_k-\bar{\sigma}_k\bar{\zeta})^2}{\bar{t}_k^2} \scalebox{0.75}[1.0]{\( - \)} \tfrac{\bar{\zeta}^2}{\bar{r}_k^2}\Big)\,\mathrm{d}\bar{\zeta} \\
& = \sqrt{\frac{\pi \,\bar{t}_k^2\,\bar{r}^2}{\bar{r}^2 + \bar{t}_k^2}}\exp\Big(\scalebox{0.75}[1.0]{\( - \)} \tfrac{\bar{z}_k^2}{\bar{r}^2+\bar{t}_k^2}\Big), \\
\end{align*}

\noindent where
\begin{align*}
\bar{\xi} &= (\xi - x_k^\mathrm{f})/\sqrt{2}\sigma_k, \\
\bar{\eta} &= (\eta - y_k^\mathrm{f})/\sqrt{2}\sigma_k,\\
\bar{\zeta} &= \zeta/\sqrt{2}\sigma_k,\\[5pt]
\intertext{and the remaining substitutions are}
T_k &= P_k/\sqrt{2}\pi^{3/2}\lambda\sigma_k,\,\,\,\,\,\,\,\,\,\,\,\,  \mathrm{[K]}\\
\bar{v}_{k, x} &= v_{k, x} \sigma_k/2\sqrt{2}\kappa,\\
\bar{v}_{k, y} &= v_{k, y} \sigma_k/2\sqrt{2}\kappa,\\
\bar{t}_k^\mathrm{f} &= \sqrt{2\kappa(t_k^\mathrm{f}-t_k^\mathrm{i})}/\sigma_k,\\
\bar{r} &= \sqrt{2\kappa(t_k^\mathrm{f} - t_k^\mathrm{i} - r)}/\sigma_k,\\[3pt]
\bar{x}_k &= (x - x_k^\mathrm{f})/\sqrt{2}\sigma_k,\\
\bar{y}_k &= (y - y_k^\mathrm{f})/\sqrt{2}\sigma_k,\\
\bar{z}_k & = z/\sqrt{2}\sigma_k,\\
\bar{t}_k & = \sqrt{2\kappa(t-t_k^\mathrm{f})}/\sigma_k.
\end{align*}

\noindent Note that $\bar{x}_k$, $\bar{y}_k$, and $\bar{t}_k$ above are slightly different from corresponding variables $\bar{x}_n$, $\bar{y}_n$, and $\bar{t}_n$ in the flux term \eqref{ustrans}. For brevity we use an index to indicate not only the segment but also if a variable is related to the flux term (index $n$) or the initial value term (index $k < n$). After some simplifications we arrive at 
\begin{equation}\label{contributions}
u_{n,k}^\mathrm{I}(\bar{\mathbf{x}}_k, \bar{t}_k) = T_k \int_0^{\bar{t}_k^\mathrm{f}} \tfrac{\bar{r}\,(\bar{r}^2 + \bar{t}_k^2)^{-1/2}}{(1 + \bar{r}^2 + \bar{t}_k^2)}  \exp \Big(\scalebox{0.75}[1.0]{\( - \)}\tfrac{(\bar{x}_k + \bar{v}_{k, x}\bar{r}^2)^2 + (\bar{y}_k + \bar{v}_{k, y}\bar{r}^2)^2}{1 + \bar{r}^2 + \bar{t}_k^2} \scalebox{0.75}[1.0]{\( - \)} \tfrac{\bar{z}_k^2}{\bar{r}^2+\bar{t}_k^2}\Big)  \,\,\mathrm{d}\bar{r}, 
\end{equation}

\noindent The above is collected into the following proposition.

\begin{proposition}\label{prop1}
Given a time partition as in Definition \ref{def1}, the solution of \eqref{pde} for piecewise constant beam parameters and piecewise linear beam path,
\[
\left\{
\begin{aligned}
\left(P(t), \sigma(t), v(t)\right) &= (P_n, \sigma_n, v_n) \\
\mathcal{C} &= \ell_n
\end{aligned}   
\right.
\,\,\,\mathrm{for}\,\, t \in \big(t_n^\mathrm{i}, t_n^\mathrm{f}\big], \,\,\, n = 1, 2, \hdots, N,\]
can be written as
\begin{align}
u(\mathbf{x}, t) &= u_n^\mathrm{I}(\mathbf{x}, t) + u_n^\Phi(\mathbf{x}, t) \notag \\   
&= u_\mathrm{init} + \sum_{k=1}^{n-1} u_{n,k}^\mathrm{I}(\mathbf{x}, t) + u_n^\Phi(\mathbf{x}, t) \,\,\, \,\,\,\mathrm{for}\,\, t \in \big(t_n^\mathrm{i}, t_n^\mathrm{f}\big], \,\,\, n = 1, 2, \hdots, N,\label{solfinal}
\end{align}

\noindent where $u_{n,k}^\mathrm{I}$ and $u_n^\Phi$ are given via \eqref{contributions} and \eqref{ustrans}, respectively.
\end{proposition}

\noindent Consequently, if $t$ lies in segment $n$, then the temperature distribution at time $t$ contains $n$ integrals over (dimensionless) time.

\section{Description of a Gauss--Legendre quadrature generation scheme}
\label{sec3}
\noindent This section is concerned with the computation of the solution \eqref{solfinal}. The solution contains a sum of integrals that need to be solved numerically. This is done using Gauss--Legendre quadrature, and particular attention is paid to the fact that the required quadrature orders may vary significantly for these integrals depending on time and beam parameters.

Gauss--Legendre quadrature requires $[-1, 1]$ to be the domain of integration. For an arbitrary function $f$ we have
\begin{align*}
\int_0^t f(s) \, \mathrm{d}s = \int_{-1}^1 \tfrac{t}{2} f\big(\tfrac{t}{2}(\tilde{s} + 1)\big)\,\mathrm{d} \tilde{s}.
\end{align*}
Such an integral can then be approximated as
\[\int_{-1}^1 \tfrac{t}{2} f\big(\tfrac{t}{2}(\tilde{s} + 1)\big)\,\mathrm{d} \tilde{s} \approx \frac{t}{2} \sum_{j = 1}^M w_j f\big(\tfrac{t}{2}(\tilde{s}_j + 1)\big),\]

\noindent where $\tilde{s}_j$ are the roots of the $M^\mathrm{th}$ Legendre polynomial $p_M(\tilde{s})$ and the corresponding weights are given by 
\[w_j = \frac{2}{(1-\tilde{s}_j^2)\big(p_M^{'}(\tilde{s}_j)\big)^2}.\]

\noindent The roots and weights can be computed easily \cite{trefethen}.

We write down the numerical approximation of \eqref{solfinal} to be used in implementation. For notational convenience, introduce symbols for the integrands in \eqref{solfinal},
\begin{align*}
f_n(\bar{\mathbf{x}}_n, \bar{s}) & =\!\begin{multlined}[t]\frac{1}{1+{\bar{s}}^2} \exp\Big(\scalebox{0.75}[1.0]{\( - \)}\tfrac{({\bar{x}_n}+\bar{v}_{n, x}{\bar{s}}^2)^2+{(\bar{y}_n}+\bar{v}_{n, y}{\bar{s}}^2)^2}{1+{\bar{s}}^2} \scalebox{0.75}[1.0]{\( - \)} \tfrac{{\bar{z}_n}^2}{{\bar{s}}^2}\Big),\\
          n = 1, 2, \hdots, N,
     \end{multlined}&\\
g_{n, k}(\bar{\mathbf{x}}_k, \bar{t}_k, \bar{r}) &= \!\begin{multlined}\frac{\bar{r}\,(\bar{r}^2 + \bar{t}_k^2)^{-1/2}}{(1 + \bar{r}^2 + \bar{t}_k^2)}  \exp \Big(\scalebox{0.75}[1.0]{\( - \)}\tfrac{(\bar{x}_k + \bar{v}_x\bar{r}^2)^2 + (\bar{y}_k + \bar{v}_y\bar{r}^2)^2}{1 + \bar{r}^2 + \bar{t}_k^2} \scalebox{0.75}[1.0]{\( - \)} \tfrac{\bar{z}_k^2}{\bar{r}^2+\bar{t}_k^2}\Big) ,\\
          k = 1, 2, \hdots, n-1.
     \end{multlined}&
\end{align*}

\noindent The approximation of the solution \eqref{solfinal} obtained using Gauss--Legendre quadrature becomes
\begin{align}
U & = U_{n}^\mathrm{I} + U_{n}^\Phi\notag\\
& = u_\mathrm{init} +  \sum_{k=1}^{n-1} T_k \bar{A}_{n, k}^\mathrm{I} + T_n \bar{A}_{n}^\Phi
\,\,\,\mathrm{for}\,\, t \in \big(t_n^\mathrm{i}, t_n^\mathrm{f}\big], \,\,\, n = 1, 2, \hdots, N,
\label{solgauss}
\end{align}
\noindent where
\begin{align*}
\bar{A}_{n, k}^\mathrm{I}(\cdot, \bar{t}_k, \bar{v}_k, \bar{t}_k^\mathrm{f})  &= \frac{\bar{t}_k^\mathrm{f}}{2}\,\sum_{j=1}^{M_{n, k}} w_{n, k, j} \,g_{n, k}\hspace{-0pt}\big(\cdot, \bar{t}_k, \tfrac{\bar{t}_k^\mathrm{f}}{2}(\tilde{r}_{j} + 1)\big), \\
\bar{A}_{n}^\Phi(\cdot, \bar{t}_n, \bar{v}_n) &= \frac{\bar{t}_n}{2}\,\sum_{j=1}^{M_n} w_{n, j} \,f_n\hspace{-0pt}\big(\cdot, \tfrac{\bar{t}_n}{2}(\tilde{s}_{j} + 1)\big).
\end{align*}

\noindent The spatial variables have been suppressed to highlight that in implementation the computation of these sums is vectorized in space. 

Note that the order of each quadrature in \eqref{solgauss} may vary. On one hand, it is important to choose quadratures of high enough order as failing to do so introduces inaccurate ripple effects in the solution. On the other hand, as the aim is to obtain approximations that all are sufficiently accurate with respect to some tolerance, 
\begin{equation}\label{tolerance}\Bigg|\int_{-1}^1 f(s)\,\mathrm{d} s - \sum_{j=1}^M w_j f(s_j)\Bigg| < \mathrm{TOL}, \end{equation} 
the required orders of the quadratures in \eqref{solgauss} will vary since the beam parameters vary and, consequently, affect the shapes of the integrands as well as the domains of integration in \eqref{solfinal}. Hence, it might be undesirable to use one global quadrature order $M$ that satisfies \eqref{tolerance} for all integrals in \eqref{solfinal} since some of those integrals would be approximated to an unnecessarily high degree. More precisely, the quadrature order $M_n$ required to approximate the integral in \eqref{ustrans} in the sense of \eqref{tolerance} depends on $\bar{t}_n$ and $\bar{v}_n$. Similarly, the quadrature order $M_{n, k}$ required to approximate the integral in \eqref{contributions} in the sense of \eqref{tolerance} depends on $\bar{t}_k$, $\bar{t}_k^\mathrm{f}$, and $\bar{v}_k$. Based on these observations, we introduce an offline step that generates two quadrature look-up tables, one for $U_{n}^\mathrm{I}$ and one for $U_{n}^\Phi$, that store the required quadrature orders. The purpose of these look-up tables is to lower the computational effort by making sure that we only use quadratures of sufficient order. Once such tables have been generated, they can be used for any type of problem in the online stage, assuming that the table sufficiently covers the parameter space. 

A simple generation of a look-up table for $U_n^\Phi$ is outlined in Algorithm 1. It uses a grid over the space $\{\bar{t}_n\} \times \{\bar{v}_n\}$ of possible parameter values. An approximation $A_n^\Phi(\cdot, \bar{t}_n, \bar{v}_n)$ of the integral in \ref{ustrans} is then computed at all grid points for increasing quadrature orders. For each grid point, once the approximation is sufficient in the sense of \eqref{tolerance} for all possible $\mathbf{\bar{x}}_n$, the quadrature order is stored in an associative array, such as Python's dictionary type, using the grid point in the parameter space as a key. This table can then be used to obtain the required quadrature order for any future approximation $\bar{A}_n^\Phi(\cdot, \bar{a}, \bar{b})$ by rounding off $(\bar{a}, \bar{b})$ to the nearest grid point and use the table. The generation of a look-up table for the terms in $U_{n}^\mathrm{I}$ is similar, the only difference being that we instead loop over a grid that covers the parameter space $\{\bar{t}^\mathrm{f}_k\} \times \{\bar{t}_k\} \times \{\bar{v}_k\}$.

There are additional ways to improve the efficiency of the solver. Since the number of integrals in $u$  increases linearly with $n$, the solution becomes more and more expensive as $n$ increases. However, the integrands in \eqref{contributions} attenuate over time as the heat diffuses, and as a result the difference between contributions from segments $k \ll n$ is small. Consequently, these segments can be combined into longer segments that are then ascribed effective beam parameters, and effective contributions can then be computed on these longer segments. This is an adaptive procedure that reduces and bounds the total number of segments and thus makes sure that the number of integrals in the solution \eqref{solfinal} do not become too big. Alternatively, if the contribution from a segment $k \ll n$ is small, the corresponding intergal could be replaced entirely by an effective constant temperature. We do not formalize these procedures here.

\vspace{0.4cm}

\begin{algorithm}[H]
 \KwData{
\newline$\bar{\mathbf{t}}$, \hfill \emph{(grid points over set of possible $\bar{t}_n$)}
\newline $\bar{\bm{v}}$, \hfill \emph{(grid points over set of possible $\bar{v}_n$)}
\newline $\mathbf{(\tilde{s}, w)}$, \hfill \emph{(set of quadratures of increasing order)}
\newline $(\tilde{s}_\mathrm{ref}, w_\mathrm{ref})$, \hfill \emph{(high order reference quadrature)}
\newline TOL. \hfill\emph{(accuracy tolerance)} 
}
\KwResult{
\begin{flushleft}\hspace{1.25cm}quadDict. \hfill \emph{(array storing quadratures for computing $\bar{A}_n^\Phi$)}\end{flushleft}
}
\Begin{
\For{$\bar{t} \in \bar{\mathbf{t}}$}{
\For{$\bar{v} \in \bar{\bm{v}}$}{
Compute $\bar{A}_n^\Phi(\cdot, \bar{t}, \bar{v}, \tilde{s}_\mathrm{ref}, w_\mathrm{ref})$\hfill\emph{(reference integral)} \\
\For{$(\tilde{s}, w) \in \mathbf{(\tilde{s}, w)}$}{
Compute $\bar{A}_n^\Phi(\cdot, \bar{t}, \bar{v}, \tilde{s}, w)$\hfill\emph{(approximated integral)} \\
$E = \,\mathrm{abs}\left(\bar{A}_n^\Phi(\cdot, \bar{t}; \bar{v}, \tilde{s}_\mathrm{ref}, w_\mathrm{ref}) - \bar{A}_n^\Phi(\cdot, \bar{t}, \bar{v}, \tilde{s}, w)\right)$\hfill\emph{(error} \\
\hfill\emph{vector)}\\
\If{$\mathrm{max}(E) < \mathrm{TOL}$}{
quadDict$(\bar{t}, \bar{v}) = (\tilde{s}, w)$\hfill\emph{(store quadrature if approxim- }\\
break\hfill\emph{ation is sufficiently accurate)}\\

}
}
}
}
}
 \caption{{\small Basic algorithm for generating a quadrature look-up table for the flux term. Since we can not compute the integral in a condition of type \eqref{tolerance} exactly, an approximation $\bar{A}_n^\Phi(\cdot, \bar{t}; \bar{v}, \tilde{s}, w)$ is instead compared to an accurate approximation based on a quadrature of very high order. The space variable has been suppressed in the algorithm to highlight that in implementation the computation of temperature is vectorized in space, and therefore we do not let the required quadrature order be space dependent. The corresponding generation for the integrals in the initial value term is similar.}}
 \label{alg1}
\end{algorithm}

\section{Numerical examples}
\label{sec4}
We use a Python implementation to look at three short examples. The first example highlights Proposition \ref{prop1}. The second example solves a larger problem and utilizes the look-up tables described in Section \ref{sec3}. Finally, the third example illustrates the importance of beam parameter optimization. Throughout this section, the initial temperature and material parameters are set to $u_\mathrm{init} = 1000 \,\mathrm{K}$, $\lambda = 20 \,\mathrm{W/mK}$, and $\kappa = 8.4495\cdot 10^{-6} \, \mathrm{m^2/s}$.

\subsection*{Example 1}
\label{ex1}
We confirm that our solution \eqref{solfinal} is correct by solving \eqref{pde} when the beam path is a single line with constant beam parameters. Let $\mathcal{C} = \{\left( x, 0, 0\right): \,\, 0 \leq x \leq 10 \, \mathrm{mm}\}$. Furthermore, let $P = 100 \,\mathrm{W}$, $\sigma = 0.1 \,\mathrm{mm}$ and $\mathbf{v} = (1000, 0) \,\mathrm{mm/s}$. The solution of this problem can be expressed differently depending on the partitioning of $\mathcal{T}$. If no partitioning is carried out, then the solution \eqref{solfinal} reduces to $u = u_\mathrm{init} + u^\Phi$. As stated in Section \ref{sec2}, the solution for this case is known from literature \cite{john} and will act as our reference solution. The reference solution is compared with another solution based on a uniform partition of $\mathcal{T}$ into four segments such that   
\[x_n^\mathrm{i} = 0.25(n-1), \,\,\, t_n^\mathrm{i} = 0.25(n-1)/v, \,\,\, (P_n, \sigma_n, v_n) = (P, \sigma, v), \,\,\,  \,\,\, n = 1, \,2,\, 3,\, 4.\]

\noindent For the sake of simpler display of results, we restrict ourselves to computing the solution on $\mathcal{C}$ only. A very high quadrature order of $800$ is used for all integrals. We use step sizes $0.05\,\mathrm{mm}$ and $0.05\,\mathrm{ms}$ in space and time, respectively. Hence, $\mathcal{C}\times\mathcal{T}$ is discretized into a grid of size $200 \times 200$. Denote the discretized reference solution by $u_\mathrm{ref}$ and the discretized solution based on the partition by $u$. The solutions align, as

\begin{equation}\sup_{\mathcal{C}\times\mathcal{T}}\left(|u_\mathrm{ref} - u|\right) = 4.4801\cdot 10^{-11} \,\,\mathrm{K}.\label{diff}\end{equation}

\noindent Recall from \eqref{contributions} the temperature contribution $u_{n,k}^\mathrm{I}$ during segment $n$ from segment $k < n$. These contributions are shown to the left in Figure~\ref{ex1} for different times. If $t$ lies in segment $n$, then $u_n^\mathrm{I} - u_\mathrm{init}$ consists of $n-1$ terms, i.e.,\ it contains one contribution from each segment completely traversed by the beam.

The total solution is shown to the right in Figure \ref{ex1} for different times. We see that the initial value term and flux term are added to form the total temperature distribution with its characteristic tail. As expected, the initial value term dominates over the flux term on the part of the beam path that does not belong to the current segment. And as shown by \eqref{diff}, the solution $u$ aligns with $u_\mathrm{ref}$. 

\newpage
\begin{figure}[H]
\centering
\begin{subfigure}{.5\textwidth}
  \centering
  \includegraphics[width=1\linewidth,trim={3cm 8cm 3cm 10cm}]{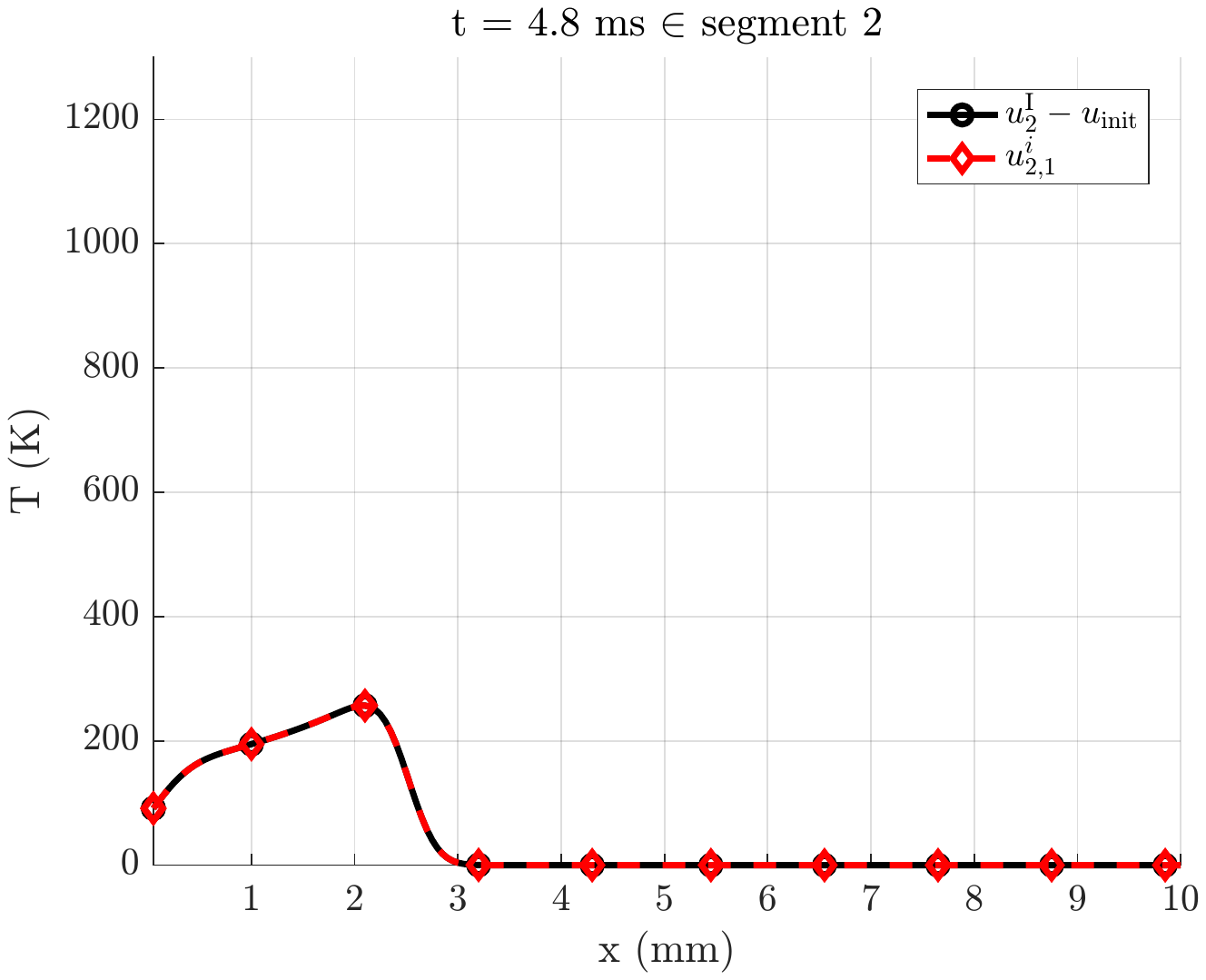}
\end{subfigure}%
\begin{subfigure}{.5\textwidth}
  \centering
  \includegraphics[width=1\linewidth,trim={3cm 8cm 3cm 10cm}]{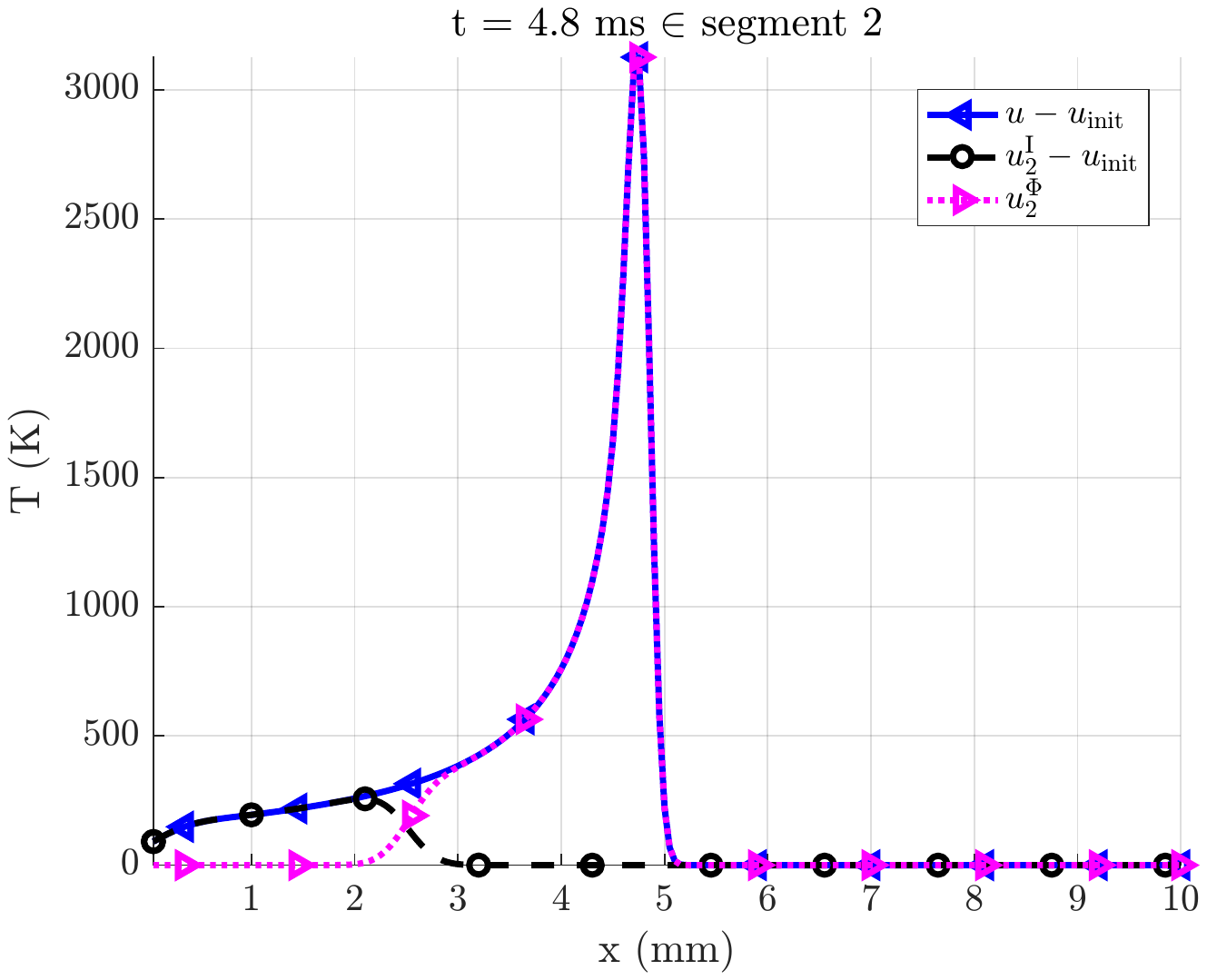}
\end{subfigure}
\begin{subfigure}{.5\textwidth}
  \centering
  \includegraphics[width=1\linewidth,trim={3cm 8cm 3cm 10cm}]{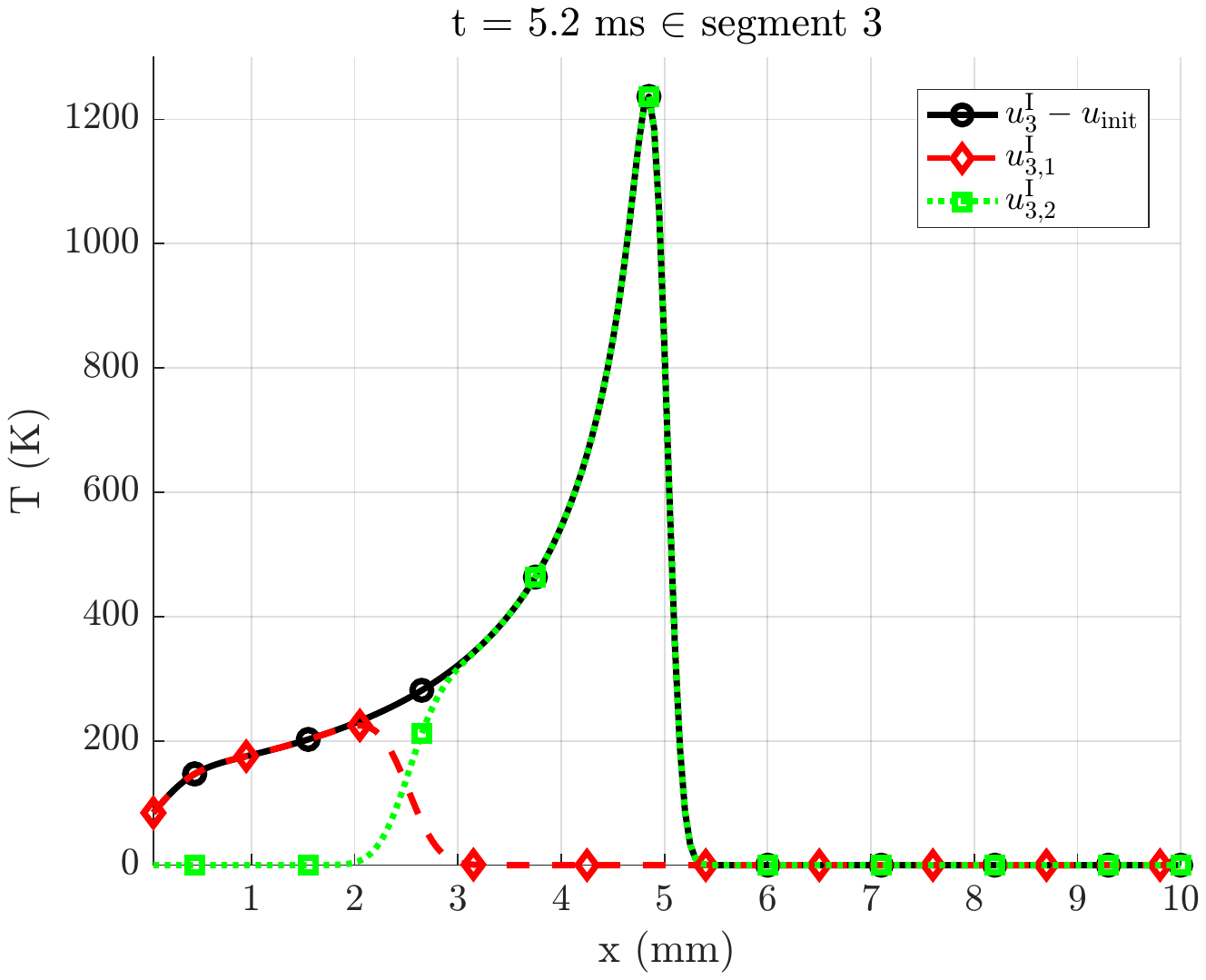}
\end{subfigure}%
\begin{subfigure}{.5\textwidth}
  \centering
  \includegraphics[width=1\linewidth,trim={3cm 8cm 3cm 10cm}]{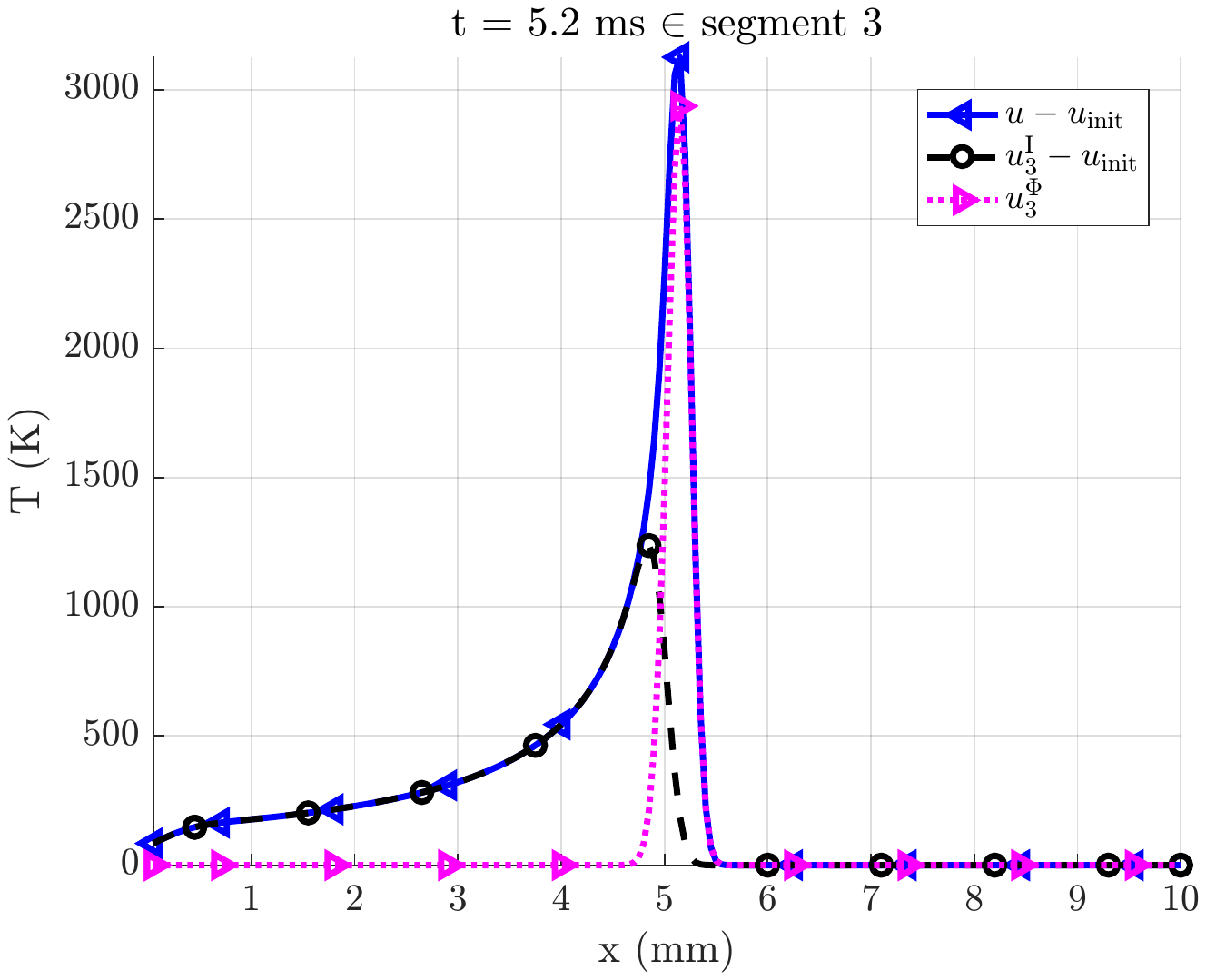}
\end{subfigure}
\begin{subfigure}{.5\textwidth}
  \centering
  \includegraphics[width=1\linewidth,trim={3cm 9cm 3cm 10cm}]{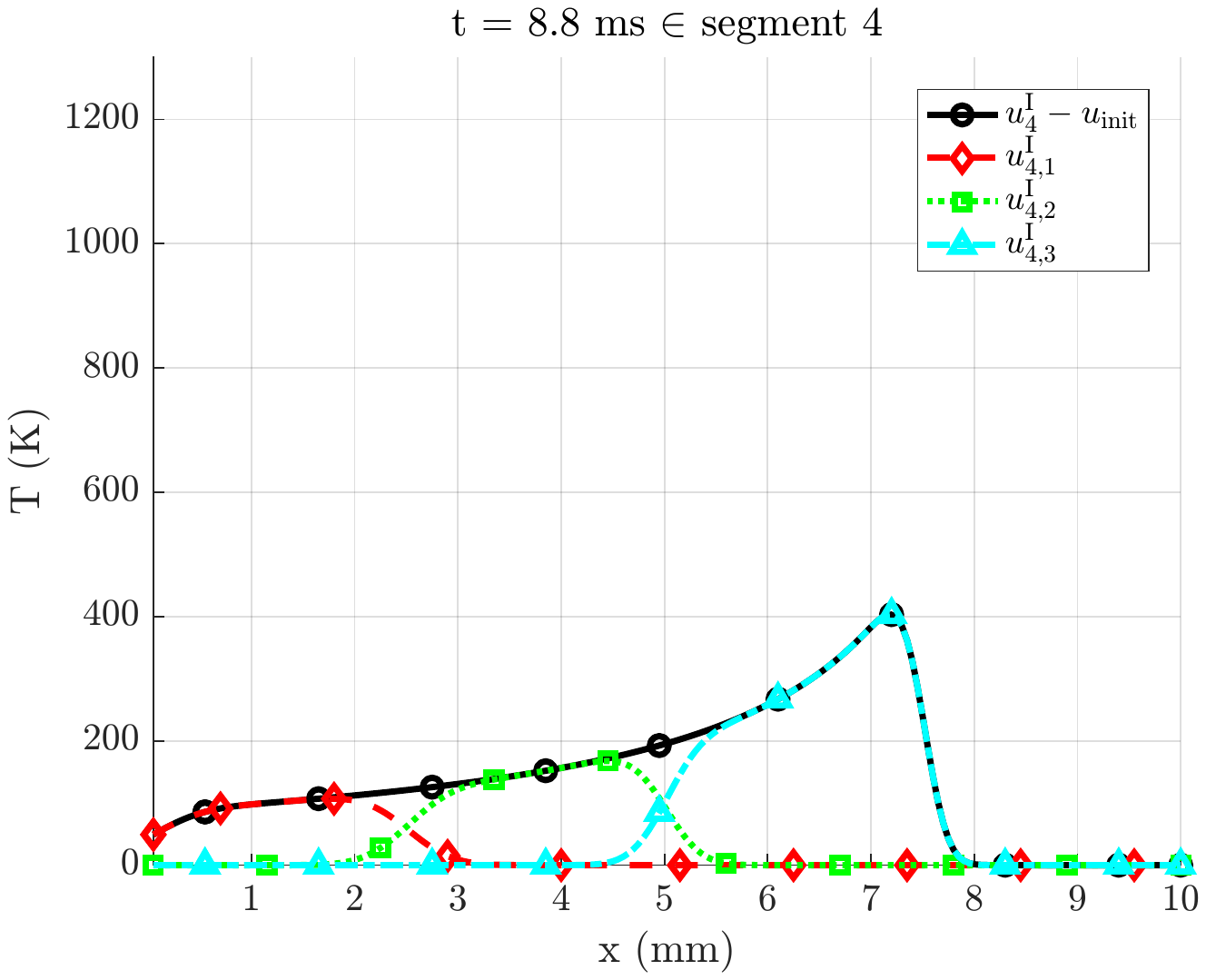}
\end{subfigure}%
\begin{subfigure}{.5\textwidth}
  \centering
  \includegraphics[width=1\linewidth,trim={3cm 9cm 3cm 10cm}]{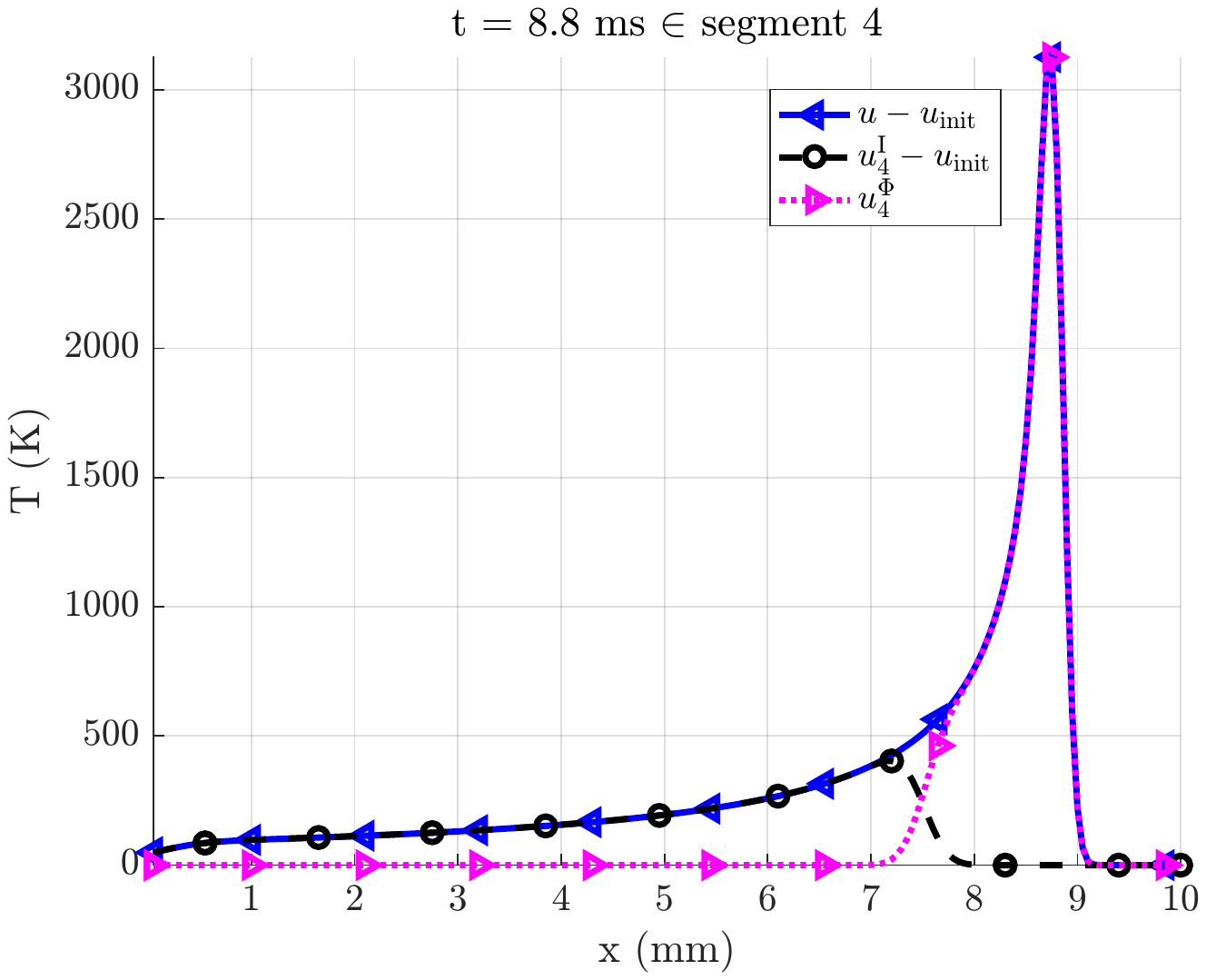}
\end{subfigure}
\caption{ \small{The solution and its components at three different times. For a time in segment $n$, the solution is given by $u(\cdot, t) = u_n^\mathrm{I}(\cdot, t) + u_n^\Phi(\cdot, t)$. The initial value term is in turn given by $u_n^\mathrm{I}(\cdot, t) = u_\mathrm{init} + \sum_{k=1}^{n-1} u_{n,k}^\mathrm{I}(\cdot, t)$.}}
\label{temp_line}
\end{figure}

\subsection*{Example 2}
\label{ex2}
\begin{figure}[H]
  \begin{minipage}[b]{0.70\linewidth}
    \centering
    \includegraphics[width=\linewidth,trim={1.5cm 14cm 1.5cm 10cm}]{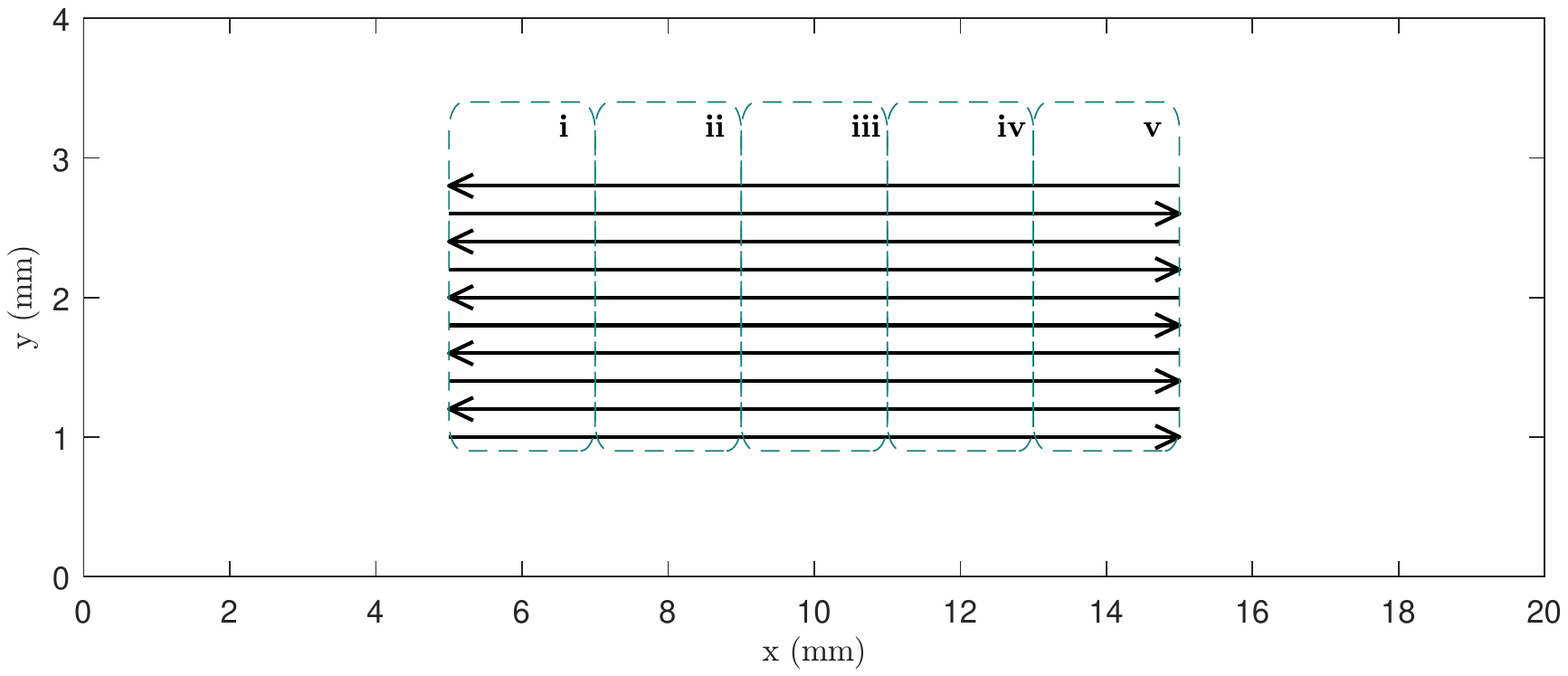}
  \end{minipage}%
  \begin{minipage}[b]{0.30\linewidth}
    \centering
    \scriptsize
\begin{tabular}{@{}c*{7}{c}c@{}}
\toprule
&\multicolumn{1}{c}{$\sigma$ (mm)}&\multicolumn{1}{c}{$v$ (m/s)}\\
\cmidrule(lr){2-2}\cmidrule(lr){3-3}
i & $0.1$ & $1$ \\
\midrule
ii & $0.15$ & $1.5$ \\
\midrule
iii & $0.2$ & $2$\\
\midrule
iv & $0.25$ & $2.5$ \\
\midrule
v & $0.3$ & $3$ \\
\bottomrule
\end{tabular}
\end{minipage}
\caption{Beam path in Example 2. It starts in the lower left corner and works its way up. The spot size and speed vary between the five sections according to the table. The effect $P = 100$ W.}
\label{beampath_snake}
\end{figure}

\begin{figure}[H]
\centering
\includegraphics[width=.78\linewidth,trim={3.5cm 11cm 1.5cm 11.5cm}]{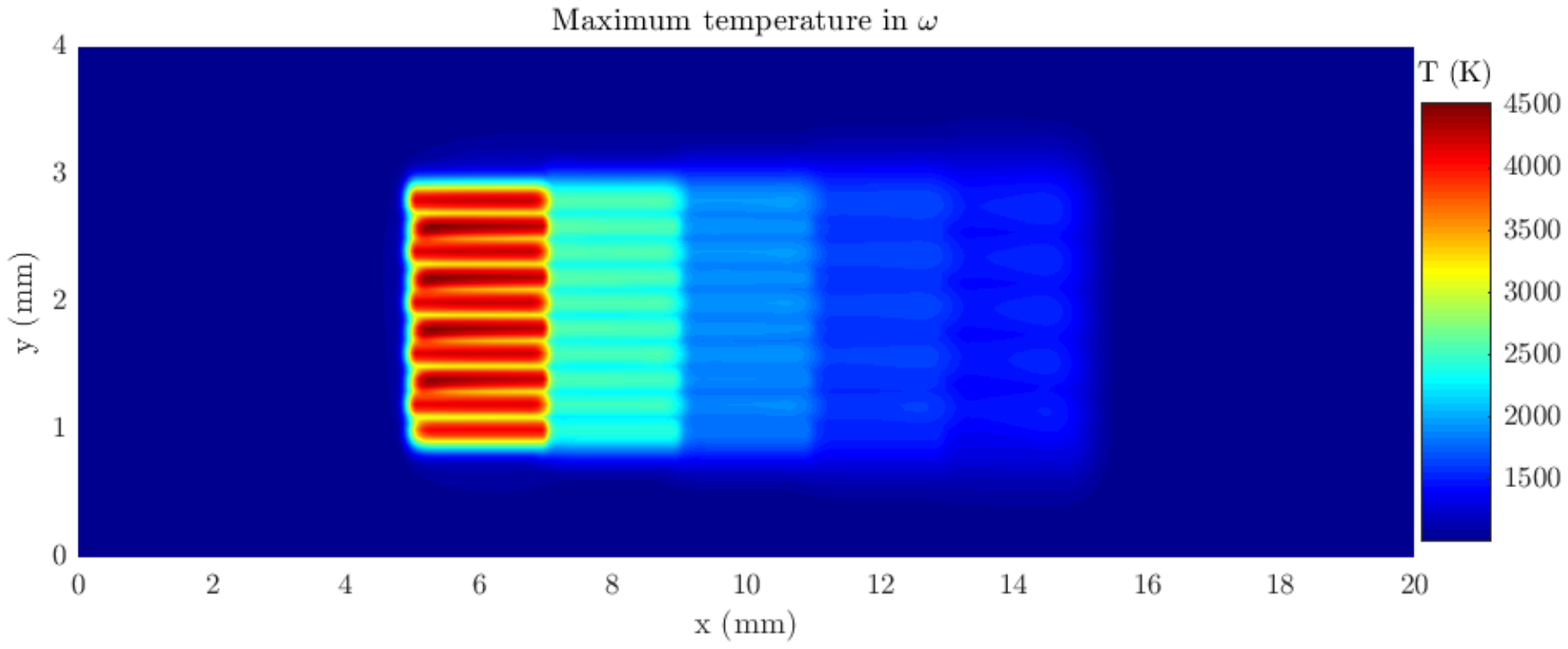}
\caption{ \small{The temperature reaches its maximum values near the turning points on the left-hand side where the beam parameters are smallest. The temperature is highly dependent on the beam parameters.}}
\label{temp_snake_max}
\end{figure}

We utilize the quadrature generation scheme described in Algorithm \ref{alg1} to compute the solution of \eqref{pde} in the subdomain $\omega = [0, 20] \times [0, 4] \times \{0\} \,  \mathrm{mm}^2$ when the beam path consists of $10$ horizontal lines of length $10 \, \mathrm{mm}$. Each line is partitioned into $5$ segments, each of length $2 \, \mathrm{mm}$, resulting in a total of $N = 50$ segments. On each of these $10$ horizontal lines, the spot size $\sigma(t)$ is increased from left to right, from $0.1\,\mathrm{mm}$ to $0.3\,\mathrm{mm}$ in steps of $0.05\,\mathrm{mm}$, and the speed $v(t)$ is increased from left to right, from $1 \, \mathrm{m/s}$ to $3 \, \mathrm{m/s}$ in steps of $0.5 \, \mathrm{m/s}$. The effect $P(t) = 100 \, \mathrm{W}$ throughout. The beam is set to traverse these $10$ lines in a snake-like fashion according to Figure \ref{beampath_snake}.
The domain $\omega$ is discretized into a square grid with a total of $401 \times 81 \times 1 = 32481$ nodes. The number of time steps is $2000$ ($40$ steps per segment). The computation of $u(\cdot, t)$ is parallelized over $t$ and the solver is executed on a laptop with an Intel Core $\mathrm{i}7$-$6500\mathrm{U}$ CPU ($4$ threads) and $8 \,\mathrm{GB}$ RAM (2133 MHz DDR4). The execution time is $683 \, \mathrm{s}$. Figure \ref{temp_snake_max} shows the maximum temperature. The temperature on the left-hand side is high since the beam spot size and speed are low there. Meanwhile, the temperature on the right-hand side is low since the beam spot size and speed are high.

\begin{figure}[H]
\centering
\begin{subfigure}{.5\textwidth}
  \centering
  \includegraphics[width=1\linewidth,trim={3cm 10cm 3cm 8cm}]{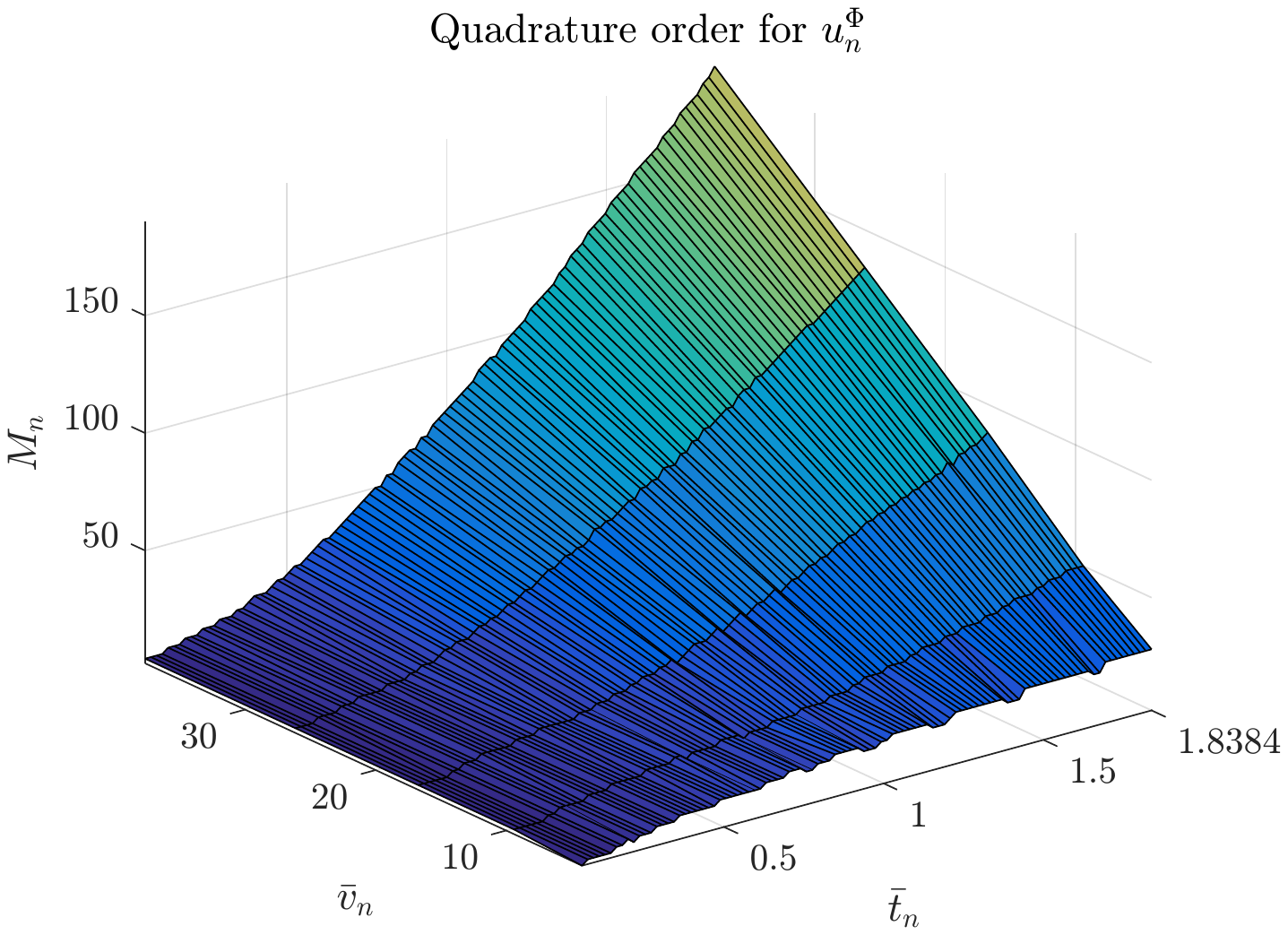}
\end{subfigure}%
\begin{subfigure}{.5\textwidth}
  \centering
  \includegraphics[width=1\linewidth,trim={3cm 10cm 3cm 8cm}]{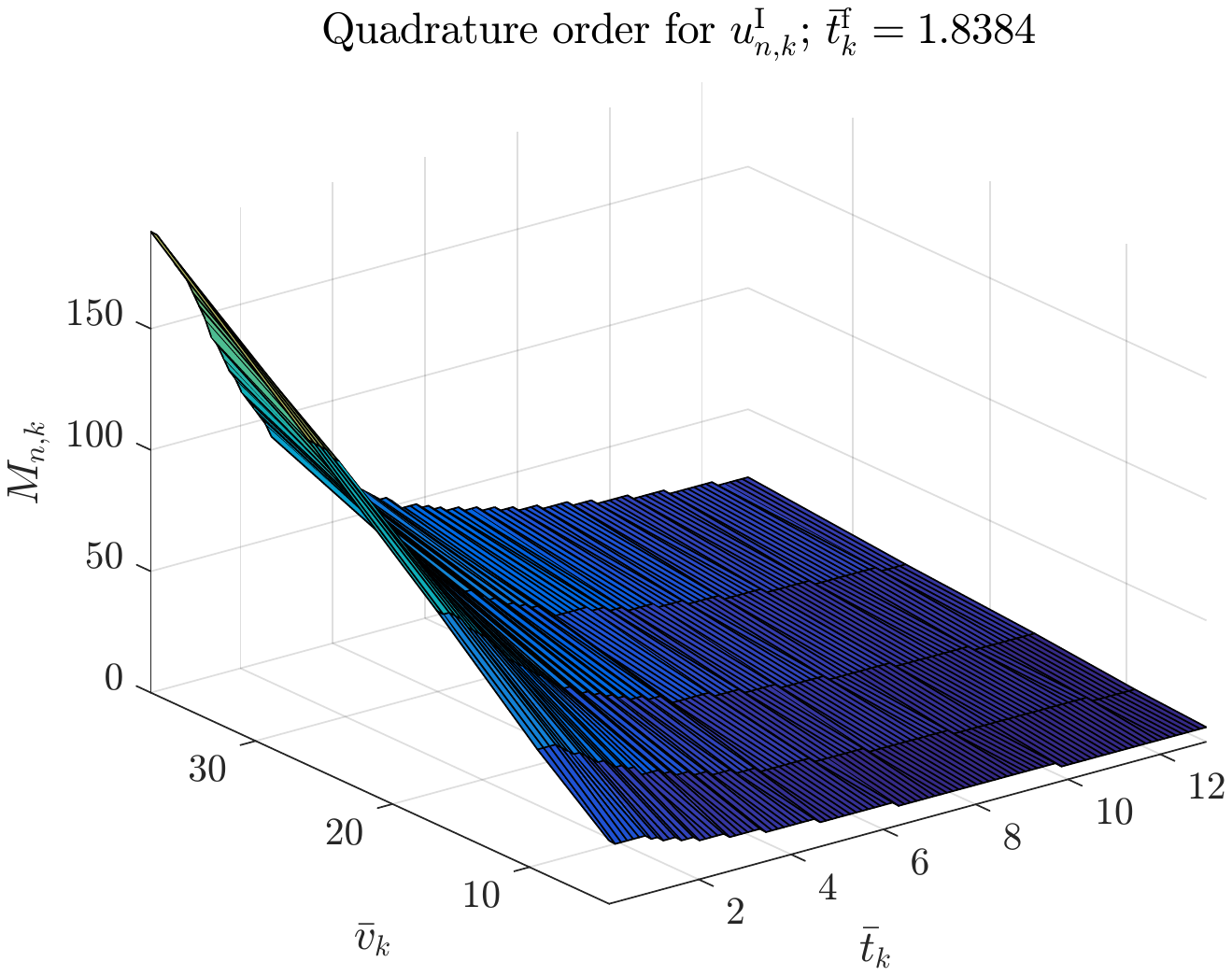}
\end{subfigure}
\caption{\small{The quadrature order required to compute $u^\Phi_n$ depends on $\bar{t}_n$ and $\bar{v}_n$. The quadrature orders required to compute the integrals $u^\mathrm{I}_{n, k}$, $k<n$, depend on $\bar{t}^\mathrm{f}_k$, $\bar{t}_k$, and $\bar{v}_k$. The required quadrature orders increase with $\bar{v}$. Furthermore, $M_n$ increases and $M_{n, k}$ decreases as their corresponding (dimensionless) time variables increase.}}
\label{quads_snake}
\end{figure}

The quadrature orders required to approximate the integrals in this example vary according to Figure \ref{quads_snake}. The approximation of all integrals are based on the tolerance $\mathrm{TOL} = 10^{-8}$ (see Algorithm \ref{alg1}).  For the integrals in the heat flux term, the required order increases with $\bar{v}_n$ and $\bar{t}_n$. We also see that as $\bar{t}_k$ increases, the required orders for the integrals in the initial value term decreases. As mentioned in Section \ref{sec3}, this decrease is due to heat diffusion and is important to utilize in big problems where the number of segments is large. For such problems, as the number of integrals in $u^\mathrm{I}$ grow, the computational cost can become detrimental unless the quadrature orders are adjusted appropriately.

\subsection*{Example 3}
\label{ex3}
One of the benefits of AM is the possibility to produce parts with complex geometries. Meanwhile, it is important to control the temperature during the melting process.  The following example illustrates the difficulty of this demand as well as the importance of beam parameter optimization.

\begin{figure}[H]
\centering
\includegraphics[width=.78\linewidth,trim={3.5cm 11cm 1.5cm 11.5cm}]{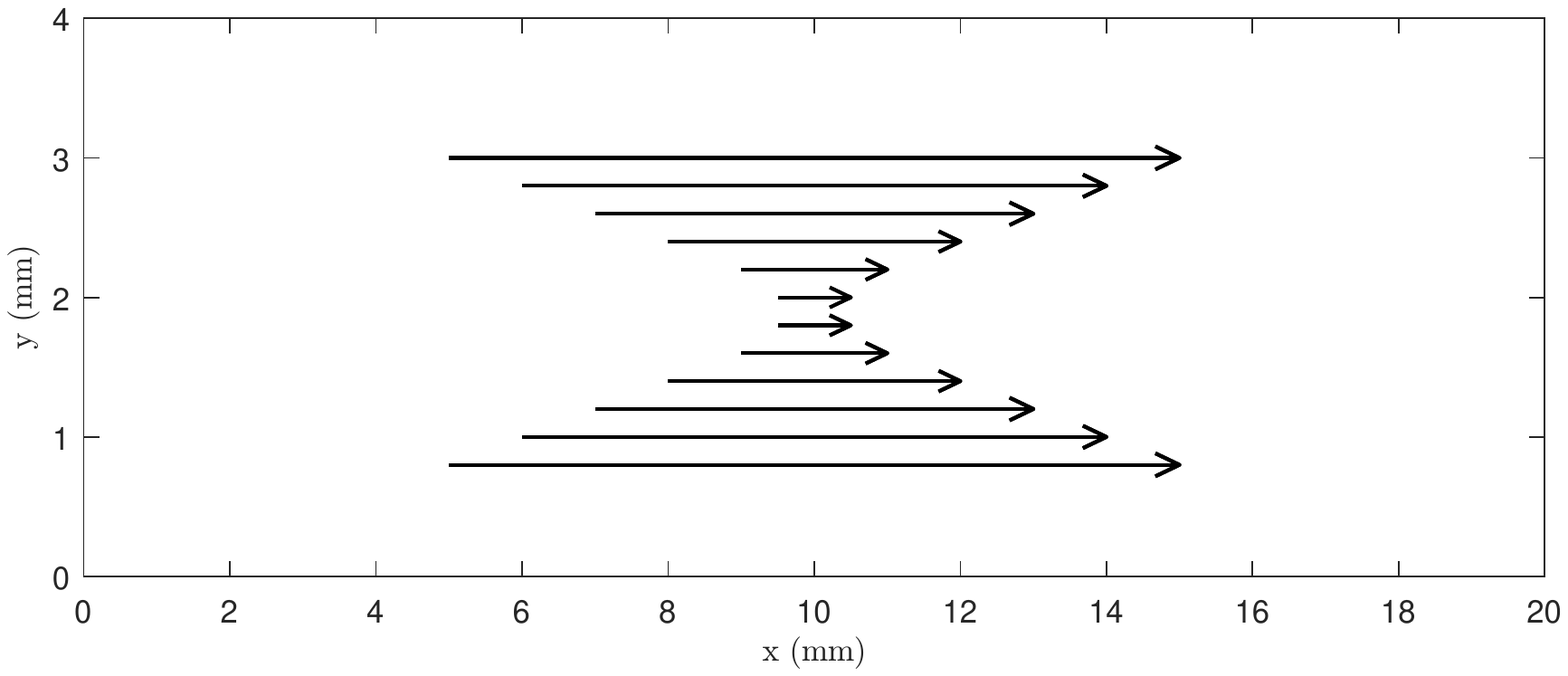}
\caption{ \small{Beam path in Example 3. It starts in the lower left corner and works its way up.}}
\label{beampath_hourglass}
\end{figure}

\begin{figure}[H]
\centering
\includegraphics[width=.78\linewidth,trim={3.5cm 11cm 1.5cm 11cm}]{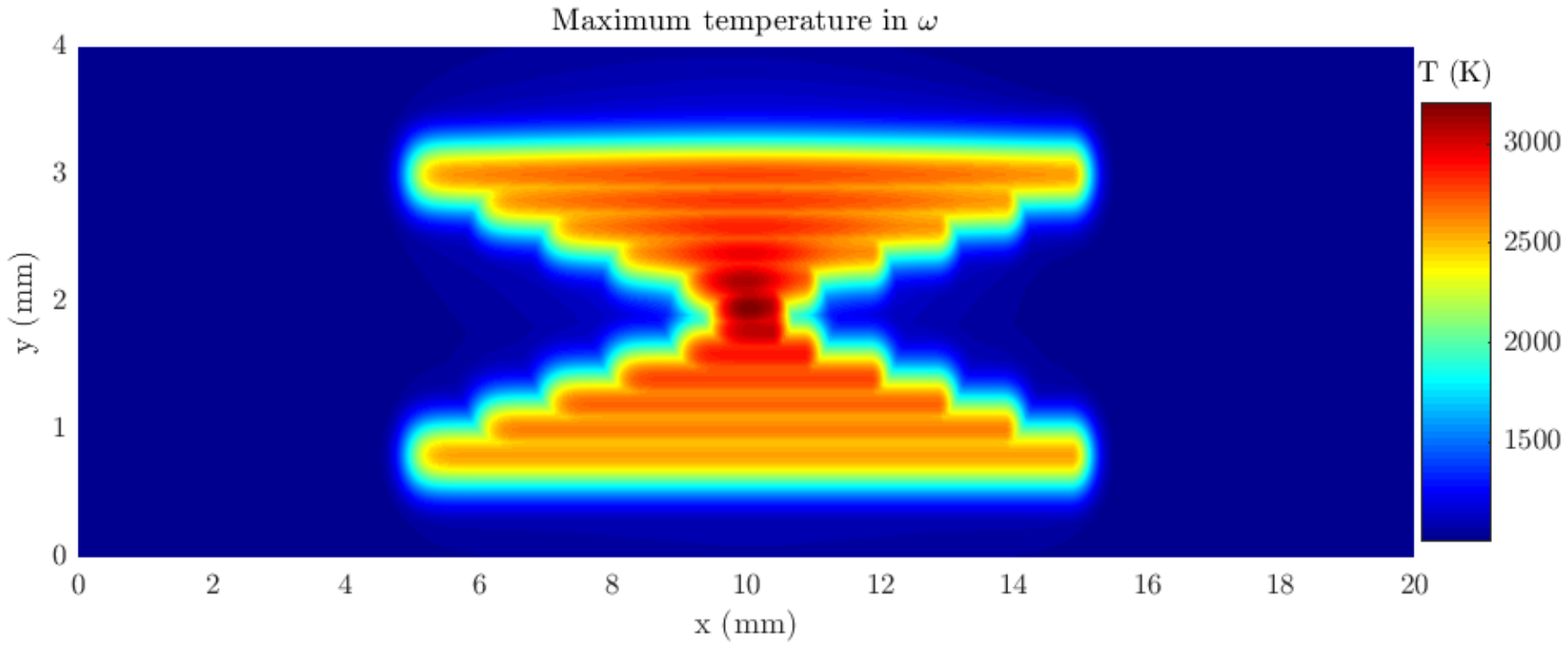}
\caption{ \small{The ratio between heat generation and heat diffusion becomes larger as the geometry of the melted area becomes smaller and more compact. Therefore, the maximum temperature increases significantly as we approach the center of the area. }}
\label{temp_hourglass_max}
\end{figure}

The beam path is shown in Figure \ref{beampath_hourglass}. The solution is evaluated on the same grid as in Example 2. We use constant beam parameters $P = 100 \,\mathrm{W}$, $\sigma = 0.1\,\mathrm{mm}$, and $v = 1000 \,\mathrm{mm/s}$. The resulting maximum temperature on $\omega$ is shown in Figure \ref{temp_hourglass_max}. The maximum temperature on the longer lines are about $2800\,\mathrm{K}$ while it exceeds $3200\,\mathrm{K}$ near the center. Now, suppose the allowed maximum surface temperature for the process is $2800 \, \mathrm{K}$. In such case, the melting displayed in Figure \ref{temp_hourglass_max} might result in too much evaporation and subsequent recoil pressure, which can result in undesired material transport such as ejection of molten materials that later cause defects \cite{gong} or formation of small ridges that prohibit the deposition of new powder layers and thus cause the manufacturing process to terminate \cite{vandenbroucke}.

We understand from Example \ref{ex2} that the high temperature in the center of the hourglass shaped melted layer can be decreased by increasing the beam spot size or speed. This amounts to an optimization procedure that is discussed further in the next, concluding section. 

\section{Conclusions and future work}
\label{sec5}
We have presented an analytical solution of the moving Gaussian heat flux problem for piecewise constant beam parameters. The intended use of the presented solution is for optimization and control of PBF processes such as electron beam melting (EBM). In EBM a large number of parameters are used to control the process. This makes it difficult to optimize the process and limits both the number of applications and the number of manufactured materials. Semi-empirical and highly efficient models that consolidate process parameters into physical expressions will provide keys to overcome the technical barriers of optimization and increase the acceptance of the technology. While our solution is efficient due to its simple form, the corresponding model does not include all the material and powder properties explicitly and can only provide an estimate of the melt pool characteristics. However, it is anticipated that effective parameters can be used to make the model reliable enough for control and optimization.    

Due to the aforementioned demand on high efficiency, we have paid particular attention to the numerical computation of the solution. Additionally, it should be noted that switching from Python to a compiled programming language such as C++ will increase the speed considerably. Adding the possibility of using highly parallel computer resources will also have a great impact on the speed of calculation. In subsequent work the intention is to do such implementations and run the calculations on high end computers. The results will be compared to both experimental melting in EBM machines and with high fidelity models and simulations.

\section*{Acknowledgments}
We acknowledge Dr. Stefan Jakobsson and Dr. Sven Johansson at Arcam EBM for their support and valuable discussions. This work was supported by the Swedish Foundation for Strategic Research under the contract ''Industrial PhD 2015 -- ID15-005''.

\newpage

\appendix
\section{Derivation of initial value term}
\label{app}
We show \eqref{ui}. In order to do this, we use that the convolution of two univariate Gaussian functions with zero means and standard deviations $\alpha_1$ and $\alpha_2$ is a Gaussian function with mean $0$ and standard deviation $\sqrt{\alpha_1^2 + \alpha_2^2}$ \cite{bromiley}. The Green's function

\begin{align*}
G(\mathbf{x}, t) &= \Big( \tfrac{1}{4 \pi \kappa t}\Big)^{3/2} \exp\Big(\scalebox{0.75}[1.0]{\( - \)}\tfrac{|\mathbf{x}|^2}{4\kappa t}\Big)\\
& = \Big( \tfrac{1}{4 \pi \kappa t}\Big)^{1/2} \exp\Big(\scalebox{0.75}[1.0]{\( - \)}\tfrac{x^2}{4\kappa t}\Big)\cdot
\Big( \tfrac{1}{4 \pi \kappa t}\Big)^{1/2} \exp\Big(\scalebox{0.75}[1.0]{\( - \)}\tfrac{y^2}{4\kappa t}\Big)\cdot
\Big( \tfrac{1}{4 \pi \kappa t}\Big)^{1/2} \exp\Big(\scalebox{0.75}[1.0]{\( - \)}\tfrac{z^2}{4\kappa t}\Big)
\end{align*}
 is a product of three univariate Gaussian functions, all with mean $0$ and standard deviation $\sqrt{2\kappa t}$. Therefore, we have
\begin{equation}\label{convrule} G(\cdot, t-t_{q}) * G(\cdot, t_{q} - t_p) = G(\cdot, t - t_p), \,\,\,\, t > t_q > t_p.\end{equation}

\begin{proposition}\label{prop2} 
Given a partition as in Definition \ref{def1}, the initial value term $u_n^\mathrm{I}(\mathbf{x}, t) = u^\mathrm{I}(\mathbf{x}, t)\big|_{t\in (t_n^\mathrm{i}, t_n^\mathrm{f}]}$, $n = 1, 2, \hdots, N$,  is given by 
\begin{equation} u_n^\mathrm{I}(\mathbf{x}, t) = u_\mathrm{init} + \sum_{k=1}^{n-1} \left(G(\cdot, t- t_k^\mathrm{f}) * u_{k}^\Phi(\cdot, t_k^\mathrm{f})\right)(\mathbf{x}). \tag{\ref{ui}}
\end{equation}
\end{proposition}
\begin{pf} In the case $n=1$, \eqref{ui} reduces to $u_1^\mathrm{i}(\mathbf{x}, t) = u_\mathrm{init}$, which is true since the initial value term on segment $1$ is equal to the initial value term in the single line solution \eqref{us_single}. Now, note that in the frame of reference of segment $n+1$, the initial temperature is the solution at the final time in segment $n$. Assuming that \eqref{ui} is true and recalling that $t_{n+1}^\mathrm{i} = t_n^\mathrm{f}$, we find 

\begin{align*}
 u_{n+1}^\mathrm{i}(\cdot, t) &\stackrel{\mathmakebox[\widthof{=}]{\eqref{u11}}}{=} \,\,\,G(\cdot, t - t_{n+1}^\mathrm{i}) * u(\cdot, t_{n+1}^\mathrm{i}) \\
 & = \,\,\,G(\cdot, t - t_{n}^\mathrm{f}) * u(\cdot, t_{n}^\mathrm{f}) \\
 &  \stackrel{\mathmakebox[\widthof{=}]{\eqref{sol}}}{=} \,\,\,G(\cdot, t - t_{n}^\mathrm{f}) * u_n^\mathrm{I}(\cdot, t_{n}^\mathrm{f}) + G(\cdot, t - t_{n}^\mathrm{f}) * u_n^\Phi(\cdot, t_{n}^\mathrm{f}) \\
 & \stackrel{\mathmakebox[\widthof{=}]{\eqref{ui}}}{=} \,\,\,\!\begin{multlined}[t]G(\cdot, t - t_{n}^\mathrm{f}) * \left(u_\mathrm{init} + \sum_{k=1}^{n-1} G(\cdot, t_{n}^\mathrm{f}- t_k^\mathrm{f}) * u_{k}^\Phi(\cdot, t_k^\mathrm{f})\right) \\ + G(\cdot, t - t_{n}^\mathrm{f}) * u_n^\Phi(\cdot, t_{n}^\mathrm{f})\end{multlined} \\
 & \stackrel{\mathmakebox[\widthof{=}]{\eqref{convrule}}}{=}\,\,\, u_\mathrm{init} + \sum_{k=1}^{n-1} G(\cdot, t-t_k^\mathrm{f}) * u_{k}^\Phi(\cdot, t_k^\mathrm{f}) + G(\cdot, t-t_{n}^\mathrm{f})*u_n^\Phi(\cdot, t_{n}^\mathrm{f})\\
 & =\,\, \,u_\mathrm{init} + \sum_{k=1}^{n} G(\cdot, t - t_k^\mathrm{f}) * u_{k}^\Phi(\cdot, t_k^\mathrm{f}).
\end{align*}
\noindent Hence \eqref{ui} holds for $n = 1, 2, \hdots, N$ by induction.$\qed$
\end{pf}

\section*{References}

\bibliography{bibliography}

\end{document}